\documentclass[final,12p,times]{elsarticle}

\usepackage{amsmath,amsthm,graphics,amssymb,amsfonts,graphicx,paralist}
\newtheorem{theorem}{Theorem}[section]
\newtheorem{corollary}{Corollary}[section]
\newtheorem{example}{Example}

\newtheorem{proposition}{Proposition}[section]

\theoremstyle{definition}

\journal{}

\begin{document}

\begin{frontmatter}



\title{On the existence and construction of Dulac functions}

\author[label1]{ Osvaldo Osuna}
\ead{osvaldo@ifm.umich.mx}
\author[label2]{ Joel Rodr\'{i}guez-Ceballos}
\ead{joel@ifm.umich.mx}
\author[label3,label4]{Cruz Vargas-De-Le\'on}
\ead{leoncruz82@yahoo.com.mx}
\author[label1]{ Gabriel Villase\~nor-Aguilar }
\ead{gabriel@ifm.umich.mx}

\address[label1]{Instituto de F\'{\i}sica y Matem\'aticas, Universidad Michoacana. Edif. C-3, Ciudad Universitaria, C.P. 58040. Morelia, Michoac\'an, M\'exico. }
\address[label2]{Instituto Tecnol\'ogico de Morelia, Departamento de Ciencias B\'asicas. Edif. AD, Morelia Michoac\'an, M\'exico
}
\address[label3]{Unidad Acad\'emica de Matem\'aticas, UAGro,
Ciudad Universitaria s/n Chilpancingo, Guerrero, M\'exico }
\address[label4]{Facultad de Ciencias, UNAM,
M\'exico, D.F. 04510, M\'exico}

\begin{abstract}
We provide sufficient conditions on the components of a vector field, which ensure the existence of Dulac
functions depending on special functions for such vector field. We also present
some applications and examples in order to illustrate our results.
\end{abstract}


\end{frontmatter}

%\begin{linenumbers}




\section{Introduction} Many problems of the qualitative theory of differential equations in the plane are related to closed orbits, this fact motivates their study. But deciding whether an arbitrary differential equation has periodic orbits or not is a difficult question that remains open.

\medskip

There are some criteria that allow us to rule out the existence of
periodic orbits in the plane; between them, we will take
particular interest in studying the Bendixson Dulac criterion. It
is well known that Bendixson-Dulac criterion is a very useful tool
for investigation of limit cycles (see \cite{Ch}, \cite{GCh},
\cite{CCGM}, \cite{GG}).

\medskip

For convenience, we recall the Bendixson-Dulac criterion, see
(\cite{F}, pag. 137).

\medskip

\begin{theorem} {\bf(Bendixson-Dulac criterion)}
Let $f_1(x_1,x_2)$, $f_2(x_1,x_2)$ and $h(x_1,x_2)$ be functions
$C^{1}$ in a simply connected domain $\Omega \subset
\mathbb{R}^{2}$ such that ${\partial (f_1h) \over \partial
x_1}+{\partial (f_2h) \over \partial x_2}$  does not change sign
in $\Omega$ and vanishes at most on a set of measure zero. Then the
system
\begin{equation}\label{syst1}
\left\{ \begin{array}{ll}
\dot{x}_1 & = f_1(x_1,x_2), \\
\dot{x}_2 & = f_2(x_1,x_2),  \; \; (x_1, \; x_2) \in \Omega,\\
\end{array} \right.
\end{equation}
does not have periodic orbits in $\Omega$.
\end{theorem}
Note that the criterion Bendixson-Dulac also discards existence of
polycycles, making it useful in establishing global stability. A
function $h$ as in the theorem is called a \textit{Dulac function}. It
is well known that Dulac function are an important tool in the
qualitative study of differential equations, but their
determination is a difficult problem. Our main goal is to give
conditions that imply the existence of Dulac functions which
depends on special forms. We give some consequences and examples to illustrate
applications of these results.

\medskip

\section{Results}

Consider the vector field
$F(x_{1},x_{2})=(f_{1}(x_{1},x_{2}),f_{2}(x_{1},x_{2}))$, then the
system \eqref{syst1} can be rewritten in the form

\begin{equation}\label{syst2}
\dot{x} = F(x),   \hspace{.5in} x=(x_{1},x_{2})\in \Omega. \\
\end{equation}
As usual the divergence of the vector field $F$ is defined by
$$
div(F)=div(f_1, f_2)=\frac{\partial f_1}{\partial
x_1}+\frac{\partial f_2}{\partial x_2}.
$$

We consider $C^{0}(\Omega, \mathbb{R})$ the set of continuous
functions and define the set

\medskip

{\small $\mathcal{F}_{\Omega} =\{f \in C^{0}(\Omega, \mathbb{R}) :
f$ doesn't change sign and vanishes only on a measure zero set\}.}

\medskip

Also for the simply connected region $\Omega$, we introduce the
sets

$$\mathcal{D}_{\Omega}^{+}(F)=\{h \in C^{1}(\Omega,\mathbb{R}) : k:= \frac{\partial (hf_{1})}{\partial x_{1}} + \frac{\partial (hf_{2})}{\partial x_{2}}\geq 0, \,\,\, k \in\mathcal{F}_{\Omega}\}$$

and

$$\mathcal{D}_{\Omega}^{-}(F)=\{h \in C^{1}(\Omega,\mathbb{R}) : k:= \frac{\partial (hf_{1})}{\partial x_{1}} + \frac{\partial (hf_{2})}{\partial x_{2}}\leq 0, \,\,\, k \in\mathcal{F}_{\Omega}\}.$$

A Dulac function in the system \eqref{syst1} of the
Bendixson-Dulac theorem is an element in the set

$$\mathcal{D}_{\Omega}(F):=\mathcal{D}_{\Omega}^{+}(F)\cup \mathcal{D}_{\Omega}^{-}(F).$$

\medskip

Our results are established with the help of the techniques
developed in \cite{OV} and
\cite{OVS}, let us recall the following result

\medskip

\begin{theorem}{\bf (\cite{OV})} If there exist $c \in \mathcal{F}_{\Omega}$ such that $h$ is a solution of the system
\begin{equation}\label{ecasociada}
f_{1}\frac{\partial h}{\partial x_{1}}+f_{2}\frac{\partial
h}{\partial x_{2}}=h\left( c(x_{1},x_{2})-div(F) \right),
\end{equation}
with $h \in \mathcal{F}_{\Omega}$, then $h$ is a Dulac function  for \eqref{syst1} on $\Omega$.
\end{theorem}

\medskip

The next theorem is our first result. In this regard Dulac
functions depending on special functions

\begin{theorem} Let $\Omega$ be a simply connected open set. Suppose a vector field
$$
F=f_1\frac{\partial}{\partial x_1}+f_2 \frac{\partial}{\partial
x_2} \in C^1(\Omega, \mathbb{R}^2).
$$
If there is $c\in \mathcal{F}_{\Omega}$ such that any of the
following conditions holds, then $\mathcal{D}_{\Omega}(F)\neq
\emptyset$:
\begin{itemize}
\item [a).-] The function
$\gamma:=\frac{c-div(F)}{f_1g_2g'_1+f_2g_1g'_2}$ depends on
$z:=g_1(x_1)g_2(x_2)$ and is continuous;

\medskip

\item [b).-] The function $\eta:=\frac{c-div(F)}{f_1g_1+f_2g_2}$
depends on $z:=k_1(x_1)+k_2(x_2)$ (with $k'_i(x_i)=g_i(x_i)$, for
$i=1, 2$) and is continuous;

\medskip

\item [c).-] The function
$\sigma:=\frac{c-div(F)}{f_1\frac{\partial z}{\partial
x_1}+f_2\frac{\partial z}{\partial x_2}}$ depends on $z:=z(x_1,
x_2)$ and is continuous.
\end{itemize}

\end{theorem}
{\bf Proof:} We consider the case $a)$, the others are analogous. We seek a Dulac function,
using the associated equation \eqref{ecasociada}.

First assume that $h$ depends only on $z:=g_1(x_1)g_2(x_2)$. Thus
the equation \eqref{ecasociada} reduces to
$$
f_1(x_1, x_2)g_2(x_2)g'_1(x_1)\frac{\partial h}{\partial
z}+f_2(x_1, x_2)g_1(x_1)g'_2(x_2)\frac{\partial h}{\partial
z}=h(z)(c(x_1, x_2)-div(F)),
$$

which is rewritten as

$$
\frac{\partial \log h}{\partial
z}=\frac{c-div(F)}{f_1g_2g'_1+f_2g_1g'_2}=\gamma(z).
$$
From our hypothesis $h:=\exp\left(\int^{z}\gamma(s)ds \right)$ is
a solution of the previous equation. Now it is easy to verify that
$h=\exp\left(\int^{z}\gamma(s)ds \right)$ is indeed a Dulac
function. The proof is complete. $\hfill \Box$

The following result is a direct consequence of theorem 2.2 and
mainly gives some particular cases

\begin{corollary} Under the conditions of theorem 2.2. If there is $c\in \mathcal{F}_{\Omega}$ such that any of the following conditions holds, then $\mathcal{D}_{\Omega}(F)\neq \emptyset$:
\item [i).-]  The function $\alpha_{i}:=\frac{c-div(F)}{f_i}$
depends only on $x_i$, for some $i \in \{1, 2\}$ and is
continuous; \item [ii).-] The function
$\beta:=\frac{c-div(F)}{x_2f_1+x_1f_2}$ depends on $z:=x_1x_2$ and
is continuous; \item [iii).-] The function
$\delta:=\frac{c-div(F)}{f_1+f_2}$ depends on $z:=x_1+x_2$ and is
continuous; \item [iv).-]  The function
$\epsilon:=\frac{c-div(F)}{c_1f_1+c_2f_2}$ depends on
$z:=c_1x_1+c_2x_2$ and is continuous; \item [v).-] The function
$\kappa:=\frac{x_{2}\left[c(x_{1},x_{2})-div(F)\right]}{f_{1}(x_{1},x_{2})-\frac{x_{1}}{x_{2}}f_{2}(x_{1},x_{2})}$
depends on $z:=\frac{x_{1}}{x_{2}}$ and is continuous.
\end{corollary}

We would like to illustrate the previous results with some
examples.

\begin{example}  We consider the classic \emph{SIS} epidemiological model
with  disease-induced death
\begin {eqnarray*}
\dot{x}_{1}&=& \lambda-\mu x_{1}-\alpha x_{2},\\
\dot{x}_{2}&=&\beta (x_{1}-x_{2})x_{2}-(\alpha+\mu+\delta) x_{2},
\end {eqnarray*}
with positive parameters, again we consider
$\mathbb{R}^2_{+}=\{(x_1, x_2)\in \mathbb{R}^2 \; :\;  x_1>0, \;
x_2>0 \}$. Denote by $ F=(f_1, f_2) $ the vector field associated
to the equation, we get
$$
div(F)=-\mu+\beta x_1-2\beta x_2-(\alpha+\mu+\delta),
$$
taking $c(x_1, x_2):=-(\mu+\beta x_2)$ we can write
$$
\alpha_2:=\frac{c-div(F)}{f_2}=\frac{-\beta x_1+\beta
x_2+(\alpha+\mu+\delta)}{(x_2)\left[\beta x_1-\beta
x_2-(\alpha+\mu+\delta)\right]}=-\frac{1}{x_2},
$$
which is continuous and depends on $z:=x_{2}$ therefore by i) of
corollary 2.1 we get $\mathcal{D}_{\mathbb{R}^2_+}(F)\neq
\emptyset$.
\end{example}

\begin{example}  Consider the Lotka-Volterra equations
\begin {eqnarray*}
\dot{x}_1&=& x_1(r_1-k_1x_1-b_{12}x_2)-h_1x_1=f_1(x_1, x_2),\\
\dot{x}_2&=& x_2(r_2-k_2x_2-b_{21}x_1)-h_2x_2=f_2(x_1, x_2),
\end {eqnarray*}
since these equations model biological systems in which two
species interact, we consider $\mathbb{R}^2_{+}=\{(x_1, x_2)\in
\mathbb{R}^2 \; :\;  x_1>0, \; x_2>0 \}$. Moreover take
$k_1k_2\geq 0$ with $k_{1}, k_{2}$ are not both zero and we
subject each population to their own proportional harvesting
effort, $h_1$ and $h_2$ respectively (i.e., the selective
harvesting). Denote by $ F=(f_1, f_2) $ the vector field
associated to the equation, we get
$$
div(F)=r_1-2k_1x_1-b_{12}x_2-h_1+r_2-2k_2x_2-b_{21}x_1-h_2,
$$
taking $c(x_1, x_2):=-(k_1x_1+k_2x_2)$ we can write
$$
\beta:=\frac{c-div(F)}{x_2f_1+x_1f_2}=\frac{-(r_1-k_1x_1-b_{12}x_2-h_1+r_2-k_2x_2-b_{21}x_1-h_2)}{(x_1x_2)\left(r_1-k_1x_1-b_{12}x_2-h_1+r_2-k_2x_2-b_{21}x_1-h_2\right)}=\frac{-1}{x_1x_2},
$$
which is continuous and depends on $z:=x_{1}x_{2}$ therefore by
ii) of corollary 2.1 we have $\mathcal{D}_{\mathbb{R}^2_+}(F)\neq
\emptyset$.

\end{example}

\begin{example}  Consider the following system of equations
\begin {eqnarray*}
\dot{x}_1&=& (\alpha_1x_1+\alpha_2x_2)(\beta_1+\beta_2 x_2^{n}+\sigma_1 x_1^{2p+1})=f_1(x_1, x_2),\\
\dot{x}_2&=& (\alpha_1x_1+\alpha_2x_2)(\beta_3+\beta_4
x_1^{m}+\sigma_2 x_2^{2q+1})=f_2(x_1, x_2).
\end {eqnarray*}
Where $\alpha_1\alpha_2> 0$ and $\sigma_1\sigma_2\geq 0$ with
$\sigma_{1}, \sigma_{2}$ are not both zero. If $n$, $m$, $p$ and
$q$ are non-negative integers. Denote by $ F=(f_1, f_2) $ the
vector field associated to the equation, we get
\begin {eqnarray*}
-div(F)&=&-\alpha_1(\beta_1+\beta_2 x_2^{n}+\sigma_1
x_1^{2p+1})-\sigma_1(2p+1)x_1^{2p}(\alpha_1x_1+\alpha_2x_2)\\&-&\alpha_2(\beta_3+\beta_4
x_1^{m}+\sigma_2
x_2^{2q+1})-\sigma_2(2q+1)x_1^{2q}(\alpha_1x_1+\alpha_2x_2),
\end {eqnarray*}
taking $c(x_1,
x_2):=(\sigma_1(2p+1)x_1^{2p}+\sigma_2(2q+1)x_1^{2q})(\alpha_1x_1+\alpha_2x_2)$
we can write
\begin {eqnarray*}
\epsilon:&=&\frac{c-div(F)}{\alpha_1f_1+\alpha_2f_2},\\&=&\frac{-\alpha_1(\beta_1+\beta_2
x_2^{n}+\sigma_1 x_1^{2p+1})-\alpha_2(\beta_3+\beta_4
x_1^{m}+\sigma_2
x_2^{2q+1})}{(\alpha_1x_1+\alpha_2x_2)\left(\alpha_1(\beta_1+\beta_2
x_2^{n}+\sigma_1 x_1^{2p+1})+\alpha_2(\beta_3+\beta_4
x_1^{m}+\sigma_2
x_2^{2q+1})\right)},\\&=&\frac{-1}{\alpha_1x_1+\alpha_2x_2},
\end {eqnarray*}
which is continuous and depends on $z:=\alpha_1x_1+\alpha_2x_2$
therefore by iv) of corollary 2.1 we have
$\mathcal{D}_{\mathbb{R}^2}(F)\neq \emptyset$.
\end{example}

\medskip

In what follows we present some consequences of Theorem 2.2 and
Corollary 2.1

\medskip
Let us consider the following system:

\begin{equation} \label{syst7}
\left \{
\begin{array}{c}
\dot{x}_1=r_1(x_1)s_1(x_2), \\
\dot{x}_2=r_2(x_1)s_2(x_2).
\end{array}
\right. 
\end{equation}
We get the next result:

\begin{corollary}
Let $\Omega \subset \mathbb{R}^2$ be a simply connected open set. If $r_2s_2' \in \mathcal{F}_{\Omega}$ and $r_1\neq 0$ (or $r_1's_1 \in \mathcal{F}_{\Omega}$ and $s_2\neq 0$) then system (\ref{syst7}) has a Dulac function on $\Omega$.
\end{corollary}
{\bf Proof:} We apply i) in Corollary 2.1  $\hfill \Box$

\medskip

Analogously we can establish the following:

\begin{corollary} Let $\Omega$ be a simply connected open set. Assume that $r_1(x_1)> 0$ ($<0$) and $s'_{2}(x_2)\in \mathcal{F}_{\Omega}$, then system
\begin{equation} \label{syst75}
\left \{
\begin{array}{c}
\dot{x}_1=r_1(x_1)r_2(x_2), \\
\dot{x}_2=s_1(x_1)+s_2(x_2),
\end{array}
\right. 
\end{equation}
has a Dulac function on $\Omega$.
\end{corollary}

\medskip

The following produce results on Dulac functions for specific systems

\medskip

\begin{proposition}
Let $g_i: \mathbb{R}^{+}\to \mathbb{R}^{+}$ be continuous functions and $a_i\in \mathbb{R}^{+}$ for $i=1, 2$, then the planar system
\begin {eqnarray*}
\dot{x}_{1}&=& x_1(a_1-x_1g_1(x_2)),\\
\dot{x}_{2}&=& x_2(a_2-x_2g_2(x_1)),
\end {eqnarray*}
supports a Dulac function on $\mathbb{R}^2_{+}:=\{(x_1, x_2)\in \mathbb{R}^2 \; :\;  x_1>0, \; x_2>0 \}$.
\end{proposition}
{\bf Proof:} First assume $a_2 \geq a_1$ and taking $c(x_1, x_2):=a_2-a_1-2x_1g_1(x_2)<0$ on $\mathbb{R}^2_{+}$ then the condition i) of corollary 2.1 is written as
$$
\alpha_{1}=\frac{c-div(F)}{f_1}=-\frac{2a_1-2x_1g_1(x_2)}{x_1(a_1-x_1g_1(x_2))}=-\frac{2}{x_1}
$$
which is continuous and depends on $z:=x_1$ therefore by i) of
corollary 2.1 we get $\mathcal{D}_{\mathbb{R}^2_+}(F)\neq \emptyset$

If $a_1\geq a_2$ we use a similar argument and get a Dulac function. $\hfill \Box$

\begin{example} The following system is a basic model of facultative
mutualism \cite{Murray}
\begin {eqnarray*}
\dot{x}_1&=& r_1x_1\left[1-\frac{x_1}{k_1+b_{12}x_2}\right]=f_1(x_1, x_2),\\
\dot{x}_2&=&
r_2x_2\left[1-\frac{x_2}{k_2+b_{21}x_1}\right]=f_2(x_1, x_2),
\end {eqnarray*}
note that satisfies the conditions of the previous proposition, therefore supports a Dulac function.
\end{example}

Similarly we can establish the following result:

\begin{proposition}
Let $g_i: \mathbb{R}^{+}\to \mathbb{R}^{+}$ be continuous functions and $a_i\in \mathbb{R}^{+}$ for $i=1, 2$,
\begin {eqnarray*}
\dot{x}_1&=&x_1\left(a_1-x_1g_1(x_2)\right)-k(x_1, x_2), \\
\dot{x}_2&=&x_2\left(a_2-x_2g_2(x_1)\right), 
\end {eqnarray*}
with $\frac{\partial k}{\partial x_{1}} \geq 0$ and $a_2\geq a_1$, then the above system supports a Dulac function on $\mathbb{R}^2_{+}:=\{(x_1, x_2)\in \mathbb{R}^2 \; :\;  x_1>0, \; x_2>0 \}$.
\end{proposition}

\begin{example} We consider the following model of two species with harvest \cite{Murray}
\begin {eqnarray*}
\dot{x}_1&=& r_1x_1\left[1-\frac{x_1}{k_1+b_{12}x_2}\right]-bx_1,\\
\dot{x}_2&=& r_2x_2\left[1-\frac{x_2}{k_2+b_{21}x_1}\right],
\end {eqnarray*}
by above proposition, the system supports a Dulac function.
\end{example}

\begin{proposition}
Let $g_i: \mathbb{R}^{+}\to \mathbb{R}^{+}$ be continuous functions and $a_i\in \mathbb{R}^{+}$ for $i=1, 2$, then the planar system
$$
\dot{x}_i=x_i\left(g_i(x_j)-a_ix_i\right), \; 1\leq i\neq j\leq 2,
$$
supports a Dulac function on $\mathbb{R}^2_{+}:=\{(x_1, x_2)\in \mathbb{R}^2 \; :\;  x_1>0, \; x_2>0 \}$.

\end{proposition}
{\bf Proof:} We verified that the conditions of Corollary 2.1, i) are satisfied, in effect
$$
-div(F)=-g_1(x_2)+2a_1x_1-g_2(x_1)+2a_2x_2,
$$
and taking $c=-(a_1x_1+a_2x_2)<0$, then
$$
\beta:=\frac{c-div(F)}{x_2f_1+x_1f_2}=-\frac{1}{x_1x_2},
$$
therefore the system supports a Dulac function. $\hfill \Box$

\begin{example} The following system is a model for mutualism \cite{Graves}
\begin {eqnarray*}
\dot{x}_1&=& x_1\left[\left(r_{1}+(r_{11}-r_{1})(1-e^{-k_{1}x_2})\right)-a_1x_1\right]=f_1(x_1, x_2),\\
\dot{x}_2&=&
x_2\left[\left(r_{2}+(r_{22}-r_{2})(1-e^{-k_{2}x_1})\right)-a_2x_2\right]=f_2(x_1,
x_2),
\end {eqnarray*}
where $r_i$, $r_{ii}$, $k_i$, $a_i\in \mathbb{R}^{+}$  are
constants and $r_{ii}>r_{i}$, $i = 1, 2$. Note that satisfies the
conditions of the previous proposition, therefore supports a Dulac
function.
\end{example}

\begin{example} Gopalsamy \cite{Gopalsamy} had proposed the following model to describe the mutualism mechanism:
\begin {eqnarray*}
\dot{x}_1&=& r_1 x_1\left[\frac{k_1+a_1x_2}{1+x_2}-x_1\right]=f_1(x_1, x_2),\\
\dot{x}_2&=& r_2
x_2\left[\frac{k_2+a_2x_1}{1+x_1}-x_2\right]=f_2(x_1, x_2),
\end {eqnarray*}
where $r_i$, $k_i$, $a_i\in \mathbb{R}^{+}$  are constants and
$a_i>k_i$, $i = 1, 2$. Depending on the nature of $k_i$ (i = 1,
2), previous system can be classified as facultative, obligate or
a combination of both.
\end{example}

%\section*{Acknowledgments} O. Osuna was partially supported by CONACYT, CIC-UMSNH. We thank the referee for various comments that helped improving the paper.

\bibliographystyle{model1a-num-names}
%%\bibliography{<your-bib-database>}



%\end{linenumbers}
\end{document}